\documentclass[a4paper,10pt]{article}
\usepackage{amsmath,amsthm,amssymb,amsfonts,epsfig,epstopdf,titling,url,array,enumerate,cite,mathrsfs,upgreek,authblk,
scrextend,mathtools,xparse}
\usepackage[bindingoffset=0.2in,left=1in,right=1in,top=1in,bottom=1in,footskip=.25in]{geometry}
\usepackage{graphicx}
\usepackage{float}
\theoremstyle{plain}
\newtheorem{thm}{Theorem}[section]
\providecommand{\keywords}[1]{\begin{addmargin}[28pt]{28pt}\noindent\textbf{Keywords:} #1 \end{addmargin}}
\newtheorem{lma}[thm]{Lemma}

\newtheorem{ppn}[thm]{Proposition}
\theoremstyle{definition}
\newtheorem{dfn}[thm]{Definition}

\newtheorem{rem}[thm]{\textit{Remark}}
\providecommand{\ams}[1]{\begin{addmargin}[28pt]{28pt}\noindent\textbf{Mathematics Subject Classification:} #1\end{addmargin}}

\title{ \textbf{Some topological results on generalized parametric metric spaces} }
\author[1]{Abhishikta Das}
\author[2]{T. Bag\footnote{corresponding author}}
\affil[1,2]{Department of Mathematics, Siksha-Bhavana,	\authorcr Visva-Bharati, Santiniketan-731235, Birbhum, West Bengal, India  
\authorcr  E-mail\textsuperscript{1}:  abhishikta.math@gmail.com	
\authorcr  E-mail\textsuperscript{2}: tarapadavb@gmail.com}

\affil[1]{Orcid Id: https://orcid.org/0000-0002-2860-424X} 
\affil[2]{Orcid Id: https://orcid.org/0000-0002-8834-7097} 

\date{}
\begin{document}
		\maketitle
\begin{abstract}
		\noindent
		
In this paper, ideas of open ball, closed ball, compact set are introduced and some related basic properties are studied. Some topological properties and some other well known results of metric spaces including Cantor's intersection theorem are established in generalized parametric metric space setting. 
 
\end{abstract}
 \ams{47H10, 54H25.} 
\section{Introduction}

It is well known that metric structure plays a pivotal role in Functional Analysis and metric fixed point theories are widely applicable in various branches of science. Due to wide application potential of metric fixed point theories attract researchers which yields a lot of generalized metric spaces(viz. 2-metric\cite{12}, b-metric\cite{10}, G-metric\cite{13}, S-metric\cite{14}, cone metric\cite{15}, parametric metric\cite{1}, parametric S-metric\cite{2}, etc.) and the study is still continued(please see \cite{3,4,5,7,8,6}, etc.). Most of the generalization have been done by modifying the `triangle inequality' of metric axioms. Following the definition of parametric metric, introduced by N. Hussain \cite{1}, in our earlier paper \cite{16} we gave a definition of generalized parametric metric. In generalized parametric metric spaces we use the general binary operation ‘o’ instead of ‘+’ and slightly change the role of ‘t’ of the definition of parametric metric. In our new approach it is possible to  achieve decomposition theorems which play crucial role to develop the results in 
generalized parametric metric spaces.  \\
In this paper, we consider the generalized parametric metric space and the notion of convergence and Cauchyness  of a sequence, completeness, bounded set, introduced by us in our earlier paper\cite{16}. Exploring the results of metric spaces and parametric spaces, we give a new idea of  open ball, closed ball, compact set, diameter of a set and establish some  well known results of metric spaces including Cantor's intersection theorem in generalized parametric metric space setting. \\
 The organization of this article is as follows. \\
 Section 2 consists of some preliminary results which we exploit in our results in this manuscript. In Section 3, we define open ball, closed ball, compact set, diameter of a set, etc. Some well known topological properties and results on completeness analogous to metric space are proved. 
\section{Preliminaries}
In this section, we recall some basic definitions and results on generalized parametric metric spaces.
\begin{dfn} \cite{16}
	Let 
$ o :[0,\infty) \times [0,\infty) \rightarrow [0,\infty) $ be a binary operation which satisfies the following conditions: 
\begin{enumerate}[(a)]
	\item $ \alpha ~ o ~ 0 = \alpha ~~ $ 
	\item $  \alpha \leq \beta \implies \alpha o \gamma \leq \beta o \gamma ~~ $  
	\item $ \alpha o \gamma =  \gamma o \alpha ~~ $ 
	\item $ \alpha o ( \beta o \gamma ) = (  \alpha o  \beta ) o \gamma ~~ $ 
\end{enumerate}	
for all $ \alpha, ~ \beta, ~ \gamma \in [0,\infty) $. 
\end{dfn}
\begin{dfn} \cite{16}
$` o $' is said to be continuous if for any sequence $ \{\alpha_n \}, ~  \{\beta_n \} $ in  $ [0,\infty) $  converging to $ \alpha $
 and $ \beta $ respectively, $  \{\alpha_n  o \beta_n \} $ converges to $ \alpha o  \beta $. 
\end{dfn}
There are some additional axioms for $ ` o$'. \\
(e)  $ ` o$' is a continuous function. \\
(f) $  \alpha < \beta $ and $ \gamma < \delta $ implies $\alpha o \gamma < \beta o \delta ,$  for all $ \alpha, ~ \beta, ~ \gamma, ~ \delta \in [0,\infty)  $.  \\
 (g) $ \alpha o \alpha > \alpha$,   for all $ \alpha \in [0,\infty) $.  
\begin{dfn}  \cite{16}
	Let $X$ be a non-empty set and $ P : X \times X \times (0,\infty) \rightarrow [0,\infty) $ be a function 
	which satisfies the following conditions: 
\begin{enumerate}[(P1)]
	\item $ (P(a,b,t)=0, ~\forall t>0)~\iff a=b$;
	\item $  P(a,b,t) = P(b,a,t) ~~\forall t>0 $ and $ \forall a,b \in X $;
	\item for $s,t>0 $ and for all $ a,b, x\in X $, $ ~ P(a,b,s+t) \leq  P(a,x,s) ~ o ~ P(b,x,t) $.
\end{enumerate}	
Then the function $P$ is said to be  generalized parametric metric  and the triple $(X,P, o)$ is called a generalized parametric metric space.
\end{dfn}
\begin{rem}  
In \cite{16} we assume the following conditions,
$$ (P4) ~~ P(a,b,t) < \alpha ~ ~ \forall t > 0, ~ \text{ for any }~  \alpha  \in (0,\infty)  ~ \text{implies} ~ a = b. $$
$$ (P5) ~~\text{for all}~ a,b \in X, ~~ P(a,b,.) ~\text{is continuous function of}~ t ~~ \forall t>0. $$  
\end{rem}
We propose the following remark on the binary operation `o'.
\begin{rem}
	If the binary operation `o' is continuous, then
	\begin{enumerate}[(i)]
		\item for any $ \alpha_1, \alpha_2 > 0 $ with $ \alpha_1 > \alpha_2 $, there exists $ \alpha_3 > 0 $ such that 
		$ \alpha_2 ~ o ~ \alpha_3 \leq \alpha_1 $. 
		\item for any $ \beta_1  > 0 $, there exists $ \beta_2, \beta_3 > 0 $ such that $ \beta_2 ~ o~ \beta_3 \leq \beta_1 $.
	\end{enumerate}
\end{rem}

\begin{ppn} \cite{16}
If $ P $ is a generalized parametric  metric on $ X $ then for all $ a, b \in X $, $ P( a,b,\cdot ) $ is non-increasing function of $ t $.
\end{ppn}
\begin{ppn} \cite{16}
	Let $( X, P, \max ) $ be a generalized parametric  metric space and $ P $ satisfies (P4). Then for each $ \alpha \in (0,\infty) $,
	\begin{equation} \label{equn 1}
		d_\alpha (a,b) = \inf \{ t > 0 : P(a,b,t) < \alpha   \} ~~ \alpha \in (0,\infty) ~ \text{and}~~ \forall a,b \in X 
	\end{equation}
is the $ \alpha $-metrics induced by the generalized parametric  metric $ P $. 
\end{ppn}
%
%
\begin{lma}\cite{16}
If $ P $ is a generalized parametric  metric on $ X $  and $ d_\alpha $ is defined as in the relation (\ref{equn 1})  then $ d_\alpha $ is non-increasing function, 
for each $ \alpha  \in (0,\infty) $.
\end{lma}
\begin{dfn} \cite{16}
Let $ ( X, P, o ) $ be a generalized parametric metric space. A sequence    $ \{ x_n \}  \subseteq X $ is said to be a 
\begin{enumerate}[(i)]
\item   convergent sequence if  there exists $ x \in X $ such that 
	$ \underset{ n \rightarrow \infty }{ lim } ~ P(x_n, x, t)= 0 ~ ~ \forall t > 0 	$. \\
	The point  $ x $ is said to be limit of $ \{ x_n \} $ and denoted by $ \underset{n\rightarrow \infty }{lim} x_n $. 
\item  Cauchy sequence if $ \underset{ m, n \rightarrow \infty }{ lim }~  P ( x_n, x_m, t ) = 0 ~ ~ \forall t> 0 $.  
\end{enumerate} 
\end{dfn}
\begin{dfn} \cite{16}
	 Let $( X, P , o)$ be a generalized parametric metric space.  $ A \subseteq X $ is said to be bounded if for each $ t > 0 $ 
	 there exist a non-negative real number $ K_t  $ such that $ P(x,y,t) \leq K_t $ for all $ x,y \in A $.
\end{dfn} 
\begin{ppn} \cite{16} \label{conv bdd}
	Every convergent sequence in a generalized parametric metric space $  ( X, P, o) $ is bounded. 
\end{ppn}
In Section 3, we define topology on generalized parametric metric spaces and  prove some topological properties.
In this aspect, some  definitions  on topological spaces  are given below.
\begin{dfn} \cite{18}
	A   topological space $ ( X, \tau ) $ is said to be 
	\begin{enumerate}[(i)]
		\item $ T_0 $ space if for any two distinct points $ x, y \in X $, there exists an open set $U$ containing $x$ such that $ y \notin U $. 
		\item $ T_1 $ space if for any two distinct points $ x, y \in X $, there exists two open sets $U$ and $V$  containing $x$ and $y$ respectively  such that $ y \notin U $ and $ x\notin V $. 
		\item Hausdorff  space if for any two distinct points $ x, y \in X $, there exists two disjoint open sets $ U $ and $V$  containing $x$ and $ y $ respectively.   
		\item regular space if for each closed set $A \subset B $ and a point $ x \in X \setminus A $, there exists two disjoint open sets $ U $ and $V$  containing $A$ and $x $ respectively.
		\item normal space if for two disjoint closed sets $A$ and $B$ in $X$, there exists two disjoint open sets $ U $ and $V$  containing $A$ and $ B $ respectively. 
	\end{enumerate}
\end{dfn}

%
\section{Main results }
Here we present the idea of open ball, closed ball, diameter of a set, completeness,  compactness  etc. and prove some other important well known theorems.

\begin{dfn}\label{dfn1}
Let $ ( X, P, o ) $ be a generalized parametric metric space. For $ t> 0$, the open ball $ B(a, \alpha, t) $ with center $ a \in X $ and radius
 $ \alpha > 0 $ is defined by
		$$ B (a, \alpha, t ) = \{ b \in X : P( a, b, t) < \alpha \}. $$ 
\end{dfn}
\begin{thm} \label{thm1}
Let $ ( X, P, o ) $ be a generalized parametric metric space. We consider the collection of subsets of $ X $, 
$$ \tau_P = \{ A \subseteq X : ~ for ~ any ~ a \in A ,~  \exists ~  \alpha > 0 ~ and ~ t>0   ~ such ~ that ~ B(a, \alpha, t) \subseteq A \}. $$   
	Then	 $\tau_P $  is a  topology on  $ X $.
			\begin{proof}
			\begin{enumerate}[(i)]
			\item Clearly, $ \phi, X \in \tau_P $.
			\item Let $ A_1, A_2 \in \tau_P $ and $ a \in A_1 \cap A_2 $. Then there exists $ t_1, t_2 > 0  $ and $ \alpha_1, \alpha_2 > 0 $ such that 
			$ B(a, \alpha_i, t_i ) \subset A_i,~ i =1, 2 $. \\
			Let $ t = \max \{  t_1, t_2 \} $ and $ \alpha = \min \{ \alpha_1, \alpha_2 \} $. Then $ B(a, \alpha, t ) \subset A_i,~ i =1, 2 $. \\ 
			 Therefore $ B(a, \alpha, t ) \subset A_1 \cap A_2 $ which implies $ A_1 \cap A_2 \in \tau_P $. 
			 \item Let $  A_i \in \tau_P, ~ i \in \Delta $ and $ a \in \cup A_i $. \\
			 Then  $ a \in A_j $,  for some $ j \in \Delta $ and there exists $ t > 0, ~ \alpha > 0 $ such that 
			 $ B(a, \alpha, t) \subset A_j \subset \cup A_i $. Therefore  	$ \cup A_i \in \tau_P $.
			 \end{enumerate}	
			\end{proof}
\end{thm}
\begin{rem}
The members of $ \tau_P $ are called open sets. 
\end{rem}
\begin{thm}\label{thm2}
Let $ P $ be a generalized parametric metric on $ X $ satisfying (P5) and `o' be continuous.  Then every open ball in $ ( X, P, o ) $  is an open set.
		\begin{proof}
		For some $ t> 0 $ and $ \alpha > 0 $ we consider the open ball $ B(a, \alpha, t)$, $ a \in X $ and let $ b \in B(a, \alpha, t) $. Then 
		$ P( a, b, t) < \alpha $. \\
		Since $ P( a, b, t) $ is non-increasing and continuous function of $t (>0) $, so there exists $ t_0 \in (0, t) $ such that 
		$ P( a, b, t_0) < \alpha $. \\
		Let $ P( a, b, t_0) = \alpha_0(>0) $. \\
		Therefore $ \alpha_0 < \alpha $ and we can find  $ \beta > 0 $ such that $ \alpha_0 < \beta < \alpha $. Again there exists $ \alpha_1 > 0 $ 
		such that $ \alpha_0 ~ o ~ \alpha_1 \leq \beta $. \\
		Now we show that $ B(b, \alpha_1, t' ) \subset B( a, \alpha, t ) $ where $ t' = (t - t_0) $.   \\
		Let $ c \in B(b, \alpha_1, t' )  $. Then $ P (b, c, t' ) < \alpha_1 $. Hence,
		$$ P(a,c,t) = P (a, c, (t-t_0) + t_0  ) \leq ~ P(a,b,t_0) ~ o ~ P( b, c, t - t_0) 
		<  ~ \alpha_0 ~ o ~ \alpha_1 
		\leq  ~ \beta < \alpha. $$
		Therefore $ c \in B(a, \alpha, t) $  and hence $ B(b, \alpha_1, t' ) \subset B( a, \alpha, t ) $ where $ t' = (t - t_0) $. 
		\end{proof}
		This completes the proof.
\end{thm}
\begin{ppn}
Let $ ( X, P, o ) $ be a generalized parametric metric space satisfying 	(P5) and `o' be continuous. Then arbitrary union of open sets is open.
\begin{proof}
	The proof is straightforward.
\end{proof}	
\end{ppn}
\begin{lma} \label{lma1}
Let $ ( X, P, o ) $ be a generalized parametric metric space and the binary operation `o' be continuous. Let  $ \alpha, \beta > 0 $ be such that 
$ \beta ~ o ~ \beta \leq \alpha $ then $ \overline{B(a, \beta, \frac{t}{2})} \subset B(a, \alpha, t) $ for each $ a \in X $ and for some $ t > 0 $.
    \begin{proof}
    For some $ t> 0 $ and $ \alpha > 0 $ we consider the open ball $ B(a, \alpha, t)$, $ a \in X $. \\
    Since `o' is continuous, so for $ \alpha > 0 $, there exists $ \beta > 0 $ such that $\beta ~ o ~ \beta \leq \alpha $. \\
    Let $ b \in \overline{B(a, \beta, \frac{t}{2})} $. \\
    Since $ B(a, \beta, \frac{t}{2}) \cap B(b, \beta, \frac{t}{2}) \neq \phi $, so there exists $ c \in  B(a, \beta, \frac{t}{2}) 
    \cap B(b, \beta, \frac{t}{2}) $. Hence,
    $$ 	P(a, b, t) \leq  ~ P(a, c,\frac{t}{2}) ~ o ~ P( b, c, \frac{t}{2} ) 
		<  ~ \beta ~ o ~ \beta 		\leq  \alpha. $$
    Therefore $ b \in  B(a, \alpha, t) $ and hence $ \overline{B(a, \beta, \frac{t}{2})} \subset B(a, \alpha, t) $.
    \end{proof}
\end{lma}
\begin{dfn}\label{dfn2}
Let $ ( X, P, o ) $ be a generalized parametric metric space. For $ t> 0$, the closed ball $ B[a, \alpha, t] $ with center $ a \in X $ and radius
$ \alpha > 0 $ is defined by
		$$ B [a, \alpha, t ] = \{ b \in X : P( a, b, t) \leq \alpha \}. $$
\end{dfn}
\begin{dfn}\label{dfn3}
Let $ ( X, P, o ) $ be a generalized parametric metric space and $ A \subset X $.
\begin{enumerate}[(i)]
\item $ A $ is said to be a closed set if for any sequence $ \{ x_n \} \subset X $ converging to $ x $ implies $ x \in A$.
\item A point $ x \in X $ is said to be a limit point of $A$ if for any $ t > 0 $, there exists $ \alpha > 0 $ such that $ ( B( x, \alpha, t ) \setminus \{ x \} ) \cap A \neq \phi $.
\item The closure of $ A $, denoted by $ \overline{A} $ is the set such that $ x \in \bar{A} $ implies there exists  a sequence in $ A $ which converges to $ x $. 
\end{enumerate} 
\end{dfn}
\begin{lma} \label{lma 2} 
Let $ A $ be a closed subset of a generalized parametric metric space $ ( X, P, o )  $ and $ x \in X $. Then $ x \notin A $ implies  $ \inf \{ P(x,a, t) : a \in A \} > 0 $ for all $ t > 0 $.
 \begin{proof}
 Suppose $ x \notin A $ such that $ \inf \{ P(x,a, t) : a \in A \} = 0 $, for some $ t_0 > 0 $. \\
 This implies for a given $ \alpha \in (0,1) $, 
 \begin{align*}
 	& \exists a' \in A ~ \text{such that } ~ P(x, a', t_0) < \alpha  \\
 	\text{i.e} ~~  & \exists a' \in A ~ \text{such that } ~ a' \in B ( x, \alpha, t_0 ).
 \end{align*}
Hence, $ B ( x, \alpha, t_0 ) \cap A $ contains a point $ a' $ of $A$ other than $ x $. Thus $ x $ is a limit point of $A$. \\
Since $ A $ is closed, so $ x \in A $. \\
This leads to a contradiction.
 \end{proof}
\end{lma}
\begin{thm}\label{thm3}
Let $ ( X, P, o )$  be a generalized parametric metric space and $ P $ satisfies (P5) and `o' be continuous.  Then every closed ball  is a closed set.
		\begin{proof}
		Let $ a \in X $. For some $ t> 0 $ and $ \alpha > 0 $, we consider the closed ball $ B[a, \alpha, t] $. \\
		Next, suppose  $ b \in \overline{B[a, \alpha, t]} $. Then $ \exists $ a sequence $\{ y_n \} $ in 	$ B[a, \alpha, t] $ such that 
		$ y_n \rightarrow  y $ as $ n \rightarrow \infty $. \\
		Again we have $ P(a, y_n, t )\leq \alpha, $ for all $ n $. \\
		Now for any $ \epsilon > 0 $, 
		$$ P(a, y, t+ \epsilon ) \leq P(a, y_n,t ) ~ o ~ P(y, y_n, \epsilon ) ~ ~ \forall n. $$
		Taking limit as $ n \rightarrow \infty $ on the both sides of the above inequality, we get 
		\begin{align*}
		&  P(a, y, t + \epsilon ) \leq  \underset{ n \rightarrow \infty } {lim}   P(a, y_n,t ) 	~ o ~   \underset{ n \rightarrow \infty } {lim} P(y, y_n, \epsilon ) \\
		\implies & P(a, y, t + \epsilon  ) \leq \alpha ~ o ~ 0 = \alpha 
     	\end{align*}
		Therefore, for any $ \epsilon > 0 $,
		\begin{equation} \label{equn closed ball}
			P(a, y, t  + \epsilon ) \leq \alpha. 
		\end{equation}  
		Let $ \{ \epsilon_k \} $ be a sequence of positive real number such that $ \epsilon_k \rightarrow 0 $ as $ k \rightarrow \infty $. \\
	So from (\ref{equn closed ball}), we can write 
	$$ P(a, y, t  + \epsilon_k ) \leq \alpha.  $$
	Since P satisfies (P5), thus taking limit as $ k \rightarrow \infty $ on the both sides of the above inequality, we get 
	$$  P(a, y, t   ) \leq \alpha. $$  
		Thus $ y \in B[a, \alpha, t] $ 	 and hence $ B[a, \alpha, t] $ is a closed set.
        \end{proof}
\end{thm}
\begin{thm}\label{thm3.5}
Let $\{ x_n \} $ and $ \{ y_n \} $ be two sequences in a generalized parametric metric space $ ( X, P, o )$ converging to $ x $ and $ y $ respectively. If $ P $ satisfies (P5) and `o' be continuous, then the sequence $ \{ P (x_n , y_n, t) $ converges to   $ P (x, y, t) ~\text{for all} ~ t > 0 $.
     \begin{proof}
      Since $\{ x_n \} $ and $ \{ y_n \} $ converges to $ x $ and $ y $ respectively, thus 
      $ \underset{ n \rightarrow \infty } {lim} P(x_n, x, t ) = 0, ~ \forall t > 0 $ and 
      $ \underset{ n \rightarrow \infty } {lim} P(y_n, y, t ) = 0, ~ \forall t > 0 $. \\
      Now choose $ \epsilon > 0 $ arbitrary. Then for all $ t > 0 $, 
      \begin{align*}
		P(x, y, t+ \epsilon ) \leq & P(x, x_n, t + \frac{\epsilon}{2}) ~ o ~  P(y, x_n, t + \frac{\epsilon}{2} ) ~ ~ \forall n  \\ 
		\leq & P(x, x_n, t + \frac{\epsilon}{2}) ~ o ~  P(y, y_n,  \frac{\epsilon}{2} ) ~ o ~ P(x_n, y_n, t  ) ~ ~ ~\forall n
	  \end{align*}
	  and 
	  \begin{align*}
	  P( x_n, y_n, t + \epsilon ) \leq & P( x_n, y, t + \frac{\epsilon}{2} )  ~ o ~  P(y, y_n,  \frac{\epsilon}{2} ) ~ ~ \forall n  \\
	  \leq & P( x_n, x,  \frac{\epsilon}{2} )  ~ o ~ P(x, y, t) ~ o ~  P(y, y_n,  \frac{\epsilon}{2} )  ~ ~~ \forall n.  
	 \end{align*}
	 Taking limit as $ n \rightarrow \infty  $ on the both sides of the above inequalities, we get
	 \begin{equation} \label{seq cont 1}
	 	P(x, y, t  + \epsilon  ) \leq \underset{ n \rightarrow \infty } {lim} P( x_n, y_n, t  ) 
	 \end{equation}  
	 and
	  \begin{equation} \label{seq cont 2}
	  	 \underset{ n \rightarrow \infty } {lim} P( x_n, y_n, t + \epsilon  )   \leq  P ( x, y, t )
	  \end{equation} 
	  for all $ t> 0 $ and for all $ \epsilon > 0 $. \\ 
	  Now we claim that the relation (\ref{seq cont 1}) implies $ \underset{ n \rightarrow \infty } {lim} P( x_n, y_n, t  )  \geq P(x, y, t ) ~~ \text{for all} ~  t > 0. $ \\
	  If not, then there exists $ t_0 > 0 $ such that   $ \underset{ n \rightarrow \infty } {lim} P( x_n, y_n, t_0  )  < P(x, y, t_0 ) $. \\
	  So there exists $ \epsilon_0 > 0 $ such that $ \underset{ n \rightarrow \infty } {lim} P( x_n, y_n, t_0  )  < P(x, y, t_0 + \epsilon_0 ) < P(x, y, t_0 ) $, since P satisfies (P5). But this contradicts the relation (\ref{seq cont 1}). \\
	  Again we claim that the relation (\ref{seq cont 2}) implies $ \underset{ n \rightarrow \infty } {lim} P( x_n, y_n, t  )   \leq  P ( x, y, t ) ~~ \text{for all} ~  t > 0. $ \\
	  If not, then there exists $ t_1 > 0 $ such that $ \underset{ n \rightarrow \infty } {lim} P( x_n, y_n, t_1  )  >  P ( x, y, t_1 )  $. \\
	  So there exists $ \epsilon_1 > 0 $ such that $ \underset{ n \rightarrow \infty } {lim} P( x_n, y_n, t_1  )  > \underset{ n \rightarrow \infty } {lim} P( x_n, y_n, t_1 + \epsilon_1 )  >  P ( x, y, t_1 )  $, since P satisfies (P5). But this is a contradiction to the relation (\ref{seq cont 2}). \\
	  Therefore we obtain $ P(x, y, t ) \leq  \underset{ n \rightarrow \infty } {lim} P( x_n, y_n, t  ) \leq  P(x, y, t )  ~~ \text{for all} ~  t > 0 $ 
	  that is \\
	  $ \underset{ n \rightarrow \infty } {lim} P(x_n, y_n, t  )  =  P ( x, y, t ) ~~ \text{for all} ~  t > 0 $.
    \end{proof} 
\end{thm}
\begin{thm}\label{thm4}
Let $ ( X, P, o )$  be a generalized parametric metric space and $ P $ satisfies (P5) and `o' be continuous.  Then $ ( X, \tau_P )  $ is a Hausdorff topological space.
         \begin{proof} 
         Let $ a, b \in X $ with $ a \neq b $. Then $ \exists t_0 > 0 $ such that $ P(a, b, t_0) \neq 0 $. \\
         Let $ P(a, b, t_0) = \alpha_0(>0) $. \\
        Then  there exists $ \alpha_1 > 0 $ such that $ \alpha_1 ~ o ~  \alpha_1 < \alpha_0 $. \\
         Now consider the open sets $ B(a, \alpha_1, \frac{t_0}{2} ) $ and  $  B( b, \alpha_1, \frac{t_0}{2} )$. For any $ z \in  
         B(a, \alpha_1, \frac{t_0}{2} ) \cap  B( b, \alpha_1, \frac{t_0}{2} )$, we have $ P(a, z, \frac{t_0}{2} ) < \alpha_1 $ and 
        $ P( b, z, \frac{t_0}{2} ) < \alpha_1 $. \\ 
         If there exists $ c \in  B(a, \alpha_1, \frac{t_0}{2} ) \cap  B( b, \alpha_1, \frac{t_0}{2} ) $, then
         $$ \alpha_0  =  ~  P(a, b, t_0)  \leq  P(a, c, \frac{t_0}{2} ) ~ o ~ P( b, c, \frac{t_0}{2} )  \leq  \alpha_1 ~ o ~ \alpha_1 < \alpha _0. $$
          This leads to a contradiction and completes the proof.
         \end{proof} 
\end{thm}
\begin{thm}\label{thm5}
Let $ ( X, P, o )$  be a generalized parametric metric space and $ P $ satisfies (P5) and `o' be continuous.  Then $ ( X, \tau_P )  $ is a $ T_0 $ space.
        \begin{proof}
        Let $ a, b \in X $ with $ a \neq b $. Then there exists $ t_0 > 0 $ such that $ P(a, b, t_0) \neq 0 $. \\
        Let $ P(a, b, t_0) = \alpha_0(>0) $. \\
        Now consider the open sets $ B(a, \alpha_0, t_0 ) $. Then clearly $ b \notin  B(a, \alpha_0, t_0 ) $.
        \end{proof} 
\end{thm}
\begin{thm}\label{thm6}
Let $ ( X, P, o )$  be a generalized parametric metric space and $ P $ satisfies (P5) and `o' be continuous.  Then $ ( X, \tau_P )  $ is a $ T_1 $ space.
        \begin{proof}
        Let $ a, b \in X $ with $ a \neq b $. Then $ \exists t_0 > 0 $ such that $ P(a, b, t_0) \neq 0 $. \\
         Let $ P(a, b, t_0) = \alpha_0(>0) $. \\
        Now consider the open balls $ B(a, \alpha_0, t_0 ) $ and $ B(b, \alpha_0, t_0 ) $. Then clearly $ b \notin  B(a, \alpha_0, t_0 ) $ and \\
        $ a \notin  B(b, \alpha_0, t_0 ) $.
        \end{proof} 
\end{thm}

\begin{thm}\label{thm7}
Let $ ( X, P, o )$  be a generalized parametric metric space and $ P $ satisfies (P5) and `o' be continuous.  Then $ ( X, \tau_P )  $ is a regular space.

        \begin{proof}
        Let $ A $ be a closed set and $ x_0 \in X \setminus A $. \\
        Then by Lemma (\ref{lma 2}),  $ \alpha (t) = \inf \{ P ( x_0, a, t ) :  a \in A  \} > 0 ~ ~ \forall t > 0 $. \\
        Let us choose a fixed $ t_0 \in (0, \infty) $. Then for $ \alpha (t_0 ) $, there exists $ \alpha_1 (t_0 ) > 0 $ such that \\ 
        $ \alpha_1 (t_0 ) ~ o ~ \alpha_1 (t_0 ) < \alpha (t_0 ) $. \\
        Now consider the open sets $ U = B(x_0, \alpha_1 (t_0 ), \frac{t_0}{2} ) $ and  $ V = \cup_{a \in A} B( a, \alpha_1 (t_0 ), \frac{t_0}{2} )$ containing $ x_0 $ 
        and $ A $ respectively. \\
         We claim that $ U \cap V = \phi $. \\
        If not, then there exists $ z \in  U \cap V $. Therefore $ P(z, a, \frac{t_0}{2} ) < \alpha_1 (t_0 )$ for some $ a_0 \in A $ 
        and  \\
        $ P( x_0, z, \frac{t_0}{2} ) < \alpha_1 (t_0 ) $. Now,
		$$  P( x_0, a_0, t_0 ) \leq  P(x_0, z,  \frac{t_0}{2} ) ~ o ~  
		P( a_0, z, \frac{t_0}{2} ) \leq  \alpha_1(t_0 ) ~ o  ~ \alpha_1(t_0 ) < \alpha(t_0 ).  $$
		This contradicts our assumption and hence  the proof is complete.
         \end{proof} 
\end{thm}
\begin{thm}\label{thm8}
Let $ ( X, P, o )$  be a generalized parametric metric space and $ P $ satisfies (P5) and `o' be continuous.  Then $ ( X, \tau_P )  $ is a normal space.
        \begin{proof}
        Let $ A $ and $ B $  be two closed sets in $X$.  \\
        Then by Lemma (\ref{lma 2}), $  \alpha(t) = \inf \{P( a, b, t ) :  a \in A, b \in B  \} > 0 ~ ~ \forall t > 0 $. \\
        Let us choose a fixed $ t_0 \in (0, \infty) $. Then for $ \alpha (t_0 ) $, there exists $ \alpha_1 (t_0 ) > 0 $ such that \\ 
        $ \alpha_1 (t_0 ) ~ o ~ \alpha_1 (t_0 ) < \alpha (t_0 ) $. \\ 
        Consider the open sets $ U = \cup_{a\in A} B(a, \alpha_1 (t_0 ), \frac{t_0}{2} ) $ and  $ V = \cup_{b \in B} B( b, \alpha_1 (t_0 ), \frac{t_0}{2} )$. Then clearly 
          $ A $ and $ B $ are contained in $U$ and $V$ respectively. \\
          We claim that $ U \cap V = \phi $. \\
        If possible suppose there exists $ z \in U \cap V $. Then for  some $ a_0 \in A $ and $ b_0 \in B $ such that  \\
        $ z \in  B(a_0, \alpha_1 (t_0 ), \frac{t_0}{2} ) $ and $ z \in B( b_0, \alpha_1 (t_0 ), \frac{t_0}{2} ) $ that is $ P(a_0, z, \frac{t_0}{2} ) < \alpha_1 (t_0 )  $ and  $ P( b_0, z, \frac{t_0}{2} ) < \alpha_1 (t_0 ) $. \\
         Therefore,
        \begin{align*}
		& P( a_0, b_0, t_0 ) \leq  P(a, z, \frac{t_0}{2} ) ~ o ~  P( b, z, \frac{t_0}{2} )  \\
		\implies & P( a_0, b_0, t_0 )  \leq   \alpha_1 (t_0 ) ~ o  ~ \alpha_1 (t_0 ) \leq \alpha (t_0 ). 
		\end{align*}
		This is a contradiction to our assumption.  \\
		Hence the proof is complete.
        \end{proof} 
\end{thm}
\begin{thm}\label{thm9}
Let $ ( X, P, o )$  be a generalized parametric metric space and $ P $ satisfies (P5) and `o' be continuous.  Then $ ( X, \tau_P )  $ is 1st countable.
    \begin{proof}
    Let $ x \in X $. We have to show that $ \beta_x = \{ B(x, \frac{1}{n}, \frac{1}{n} ) : n \in \mathbb{N} \} $ is a local base for $ x\in X $. \\
    Let $ G \in \tau_P  $ such that $ x \in G $. Since $ G $ is open, then there exists $ \alpha_0 > 0 $ and $ t_0 > 0 $ such that $ B ( x, \alpha_0, t_0 ) 
    \subset G $. \\
    Choose $ n \in \mathbb{N} $ such that $ \frac{1}{n} < \alpha_0 $ and $ \frac{1}{n} < t_0 $. \\
    Then for  $ z \in X $, 
    $ z \in B(x, \frac{1}{n}, \frac{1}{n} )  $ implies $ P(x, z, \frac{1}{n} ) < \frac{1}{n} < \alpha_0 $.  Hence,
    \begin{align*}
     & P(x,z, t_0) \leq  P(x,z, \frac{1}{n} ) <  \frac{1}{n} ~~ \text{(since  P  is  non-incrasing)} \\
      \implies & P(x,z, t_0) < \alpha_0 \\
      \implies & z \in B ( x, \alpha_0, t_0 ). 
    \end{align*}
    Hence, $ B(x, \frac{1}{n}, \frac{1}{n} ) \subset B ( x, \alpha_0, t_0 ) \subset G $. Consequently $ \beta_x $ is countable local base for $ x $. \\
     Therefore  $ ( X, \tau_P )  $ is 1st countable.
    \end{proof} 
\end{thm}
\begin{thm}\label{thm10}
Let $ ( X, P, o )$  be a generalized parametric metric space and $ P $ satisfies (P5) and `o' be continuous.  Then $ ( X, \tau_P )  $ is 2nd countable.
    \begin{proof}
    Let $ \mathscr{A} = \{ a_n : n \in \mathbb{N} \} $ be a countable dense subset of $ X $ and  $ \beta = \{ B( a_m, \frac{1}{n}, \frac{1}{n} )
     : m, n \in \mathbb{N} \} $. \\
     We claim that $ \beta $ is a countable base of $ \tau_P $. \\
     Clearly $ \beta $ is countable. \\
     Let $ G $ be an open set in $ X $. Then for any $ x \in G $, there exists $ \alpha > 0 $ and $ t_0 > 0 $ such that $ B ( x, \alpha, t_0) 
     \subset G $. \\
     Again there exists $ \delta > 0 $ such that $ \delta ~ o ~ \delta < \alpha $. Let $ m \in \mathbb{N} $ such that $  \frac{1}{m} < \delta $ 
     and $ \frac{1}{m} < \frac{t_0}{2} $. \\
     Since $ \mathscr{A} $ is dense in $ X $, there exists $ a_j \in  \mathscr{A} $ such that $ a_j \in B ( x,  \frac{1}{m},  \frac{1}{m} ) $. \\
     Now for  $ a \in X $, if $ a \in  B ( x,  \frac{1}{m},  \frac{1}{m} ) $, then
    $$ P(a, x, t_0) \leq  P(a, a_j, \frac{t_0}{2} ) ~ o ~ P(x,  a_j, \frac{t_0}{2} ) 
     \leq \frac{1}{m} ~ o ~ \frac{1}{m} 
     <  \delta ~ o ~ \delta < \alpha. $$
     Hence $ a \in B(x, \alpha, t_0 ) $, and $ \beta $ is base for $ \tau_P $.
    \end{proof} 
\end{thm}
\begin{dfn}\label{dfn4}
In a generalized parametric metric space,  a set $B $ is said to be compact if  every sequence in $B$  has a  subsequence which converges  in $B$. 
\end{dfn}
\begin{thm}\label{thm11}
Let $ ( X, P, o )$  be a generalized parametric metric space and $ P $ satisfying (P5) and `o' be continuous.  Then every compact subset is closed and bounded.
    \begin{proof}
    Let  $\{x_n\} $ be a sequence in $ A   $ converging to some $ x \in X $ .  \\
    Since $ A $ is a compact set in $X$, so there exists a subsequence $ \{x_{n_k } \} $ of $ \{ x_n \}$ which converges to a point in $ A $. But $ \{x_{n_k } \} $  must converges to $x$. Hence $ x \in A $.   \\
    Since $\{x_n\} $  is an arbitrary sequence in $ A $, thus $ A $ is closed.\\
    Next, if possible suppose that $A$ is unbounded.\\
    Choose $ x_0 \in A $ any fixed element. \\
    Since $ A $ is unbounded, there exist $ x_1 \in A $ such that $ P ( x_1, x_0, t_0 ) > 1 $ for some $ t_0 > 0 $.  Similarly, there exists $ x_2 \in A $ such that $P ( x_2, x_0, t_0 ) > 2 $.
    Continuing in this way, there exists $ x_n \in A $ such that 
    $ P ( x_n, x_0, t_0  ) > n $ for all $ n \in \mathbb{N} $. \\ 
    SO we obtain a sequence $ \{ x_n \} \subseteq A $ satisfying $ P ( x_n, x_0, t_0  ) > n $ for all $ n \in \mathbb{N} $.  Since $ A $ is compact, so there exists a subsequence  $\{ x_{n_k}\}$ of $\{ x_n\}$ which convergence to some element $ x $(say) in $ A $. \\
    Therefore,  $ \underset{n\rightarrow\infty }{lim} P( x_{n_k}, x, t) = 0, ~ ~ \forall t > 0 $. \\ 
    In particular, $ \underset{n\rightarrow\infty }{lim} P( x_{n_k}, x, \frac{t_0}{2}) = 0 $. \\
    Again we have, $ P ( x_{n_k}, x_0, t_0 ) > n_k $. Now,
    $$ P ( x_{n_k}, x_0,  t_0 ) \leq P( x_{n_k}, x, \frac{t_0}{2}) ~ o ~ P( x_0, x, \frac{t_0}{2}) $$
    gives 
    $$  n_k < P( x_{n_k}, x, \frac{t_0}{2}) ~ o ~ P( x_0, x, \frac{t_0}{2}). $$
    Taking limit as $ k \rightarrow \infty $  on both sides of the above inequality, we obtain 
    $$ \infty \leq  P( x_0, x, \frac{t_0}{2}). $$
    This  contradicts that $ P $ is a real valued function.\\
    Hence $ A $ is bounded.
    \end{proof} 
\end{thm}
\begin{dfn}\label{dfn5}
Let $ ( X, P, o ) $ be a  generalized parametric metric space and $ A \subseteq X $.  Diameter of $A$, denoted by $ \delta(A) $ and defined by 
$ \delta(A) = \underset{x, y \in A }{Sup} ~ \underset{t>0 }{Sup} ~ P(x, y, t) $.
\end{dfn}
\begin{lma} \label{thm12}
Let $ ( X, P, o )$  be a generalized parametric metric space and $ P $ satisfies (P5) and `o' be continuous. For any set $ A $ in $ X $, 
 $ \delta(A) = \delta(\overline{A}) $. 
        \begin{proof}
        Since $ A \subset \overline{A} $, so 
        \begin{equation}\label{equn diam 1}
         \delta(A) \leq \delta(\overline{A}). 
        \end{equation}
        Let $ \delta(A) < \alpha_0 $. Then
        $$ \underset{x, y \in A }{Sup} ~ \underset{t>0 }{Sup} ~ P(x, y, t) < \alpha_0 $$
        \begin{equation}\label{equn diam 2}
        \implies  \underset{t>0 }{Sup} ~ P(x, y, t) < \alpha_0, ~ ~~ \forall x, y \in A 
        \end{equation}
        Choose $ x_0, y_0 \in \overline{A} $. Then $ \exists $ a sequence $ \{ x_n \} $, $ \{ y_n \} \subset A $  such that 
        $ x_n \rightarrow x_0 $ and $ y_n \rightarrow y_0 $ as $ n \rightarrow \infty $. Hence we have,
         \begin{equation}\label{equn diam 3}
         \underset{n \rightarrow \infty }{lim} P(x_n, y_n, t) = P(x_0, y_0, t ) ~ ~~ \forall t > 0 
         \end{equation}
         From (\ref{equn diam 2}) we have, 
          \begin{align*}
         & \underset{t>0 }{Sup} ~ P(x_n, y_n, t) < \alpha_0 ~ ~ ~ ~\forall n \\
          \implies & P(x_n, y_n, t) < \alpha_0 ~~ ~ ~ \forall n,~ ~\forall t > 0 \\
           \implies & \underset{n \rightarrow \infty }{lim} P(x_n, y_n, t) \leq  \alpha_0 ~ ~~ ~\forall t > 0 \\
           \implies & P(x_0, y_0, t ) \leq  \alpha_0 ~~~ ~ \forall t > 0 \\
           \implies &  \underset{t>0 }{Sup} ~ P(x_0, y_0, t ) \leq  \alpha_0 
           \end{align*}
           Since $  x_0, y_0 \in \overline{A} $ arbitrarily chosen, so  $ \underset{x, y \in \overline{A} }{Sup} ~ \underset{t>0 }{Sup} ~ P(x, y, t) \leq \alpha_0 $  i.e $ \delta(\overline{A})  \leq  \alpha_0 $ and this implies
            \begin{equation}\label{equn diam 4}
            \delta(\overline{A})  \leq  \delta(A). 
            \end{equation}
            From the relations (\ref{equn diam 1}) and (\ref{equn diam 4}), we have $ \delta(A) = \delta(\overline{A}) $.
        \end{proof} 
\end{lma}
\begin{thm}  \label{thm13}
(Cantor's intersection type theorem) Let $ ( X, P, o )$  be a generalized parametric metric space and $ P $ satisfying (P5) and `o' be continuous.
A necessary and sufficient condition that the generalized parametric metric space $ (X, P, o) $ be 
	complete is that every nested sequence of non-empty closed subsets $ \{ A_i \} $ with $ \delta(A_i) \rightarrow 0 $ as 
	$ i \rightarrow \infty $ be such that $ A = \cap _{i=1} ^\infty A_i $ contains exactly one point.
     \begin{proof}
     	Using the Lemma \ref{thm12}, the proof follows as  similar manner as the Cantor's intersection  theorem in metric spaces\cite{17}.
     \end{proof} 
\end{thm}
\begin{thm}\label{thm14}
Let $ ( X, P, o )$  be a generalized parametric metric space and $ P $ satisfying (P5) and `o' be continuous. If every Cauchy sequence in $X$ has a convergent subsequence then $X $ is complete.
     \begin{proof}
     Let $ \{ x_n \} $ be a Cauchy sequence in $ X $ and $ \{ x_{i_n}\} $ be a subsequence of $ \{ x_n \} $ which converges to $  x \in X $. \\
     We now show that $ \{ x_n \} $ also converges to $x$. \\
     Now, the inequality gives 
     $$ P(x_n, x, t) \leq  P(x_n, x_{i_n}, \frac{t}{2} ) ~ o ~ P( x_{i_n}, x, \frac{t}{2} ) ~~  ~ \forall t > 0 $$
     Taking limit as $ n \rightarrow \infty $, 
     \begin{align*}
         & \underset{n \rightarrow \infty }{lim} P(x_n, x, t) \leq 0 ~~ ~ \forall t > 0 \\
           \implies & \underset{n \rightarrow \infty }{lim} P(x_n, x, t) = 0 ~~ ~ \forall t > 0 
           \end{align*}
           Therefore the Cauchy sequence $ \{ x_n \} $ converges to $x \in X $. \\
           This completes the proof.
     \end{proof} 
\end{thm}
We end our discussion with the following result that shows the generalized parametric topology and the induced $ \alpha$-metric topology are identical.   
\begin{thm}\label{thm18}
Let $ ( X, P, \max )$  be a generalized parametric metric space such that P satisfies (P4) and $ \{ d_\alpha : ~ \alpha \in ( 0, \infty ) \} $ be the induced $ \alpha$-metric  on $ X $. Then the topology $ \tau_{d_\alpha } $ induced by $ d_\alpha $  and the topology $ \tau_P $ induced by $ P $ are identical.
        \begin{proof}
        We know for all $ a, b \in X $,
        $$ d_\alpha ( a, b ) = \inf \{ t > 0 : P(a, b, t ) < \alpha \}, ~  \alpha \in (0, \infty) $$
        Let $ G \in \tau_{d_\alpha } $. Then for each $ x \in G, $ there exists $ \epsilon > 0 $ such that  $ B_\alpha (x, \epsilon ) \subset G $ where \\
         $ B_\alpha (x, \epsilon ) = \{ y\in  X : d_\alpha (x, y ) < \epsilon \} $.  \\
         Next, let $ y \in B_\alpha (x, \epsilon ) $. Then, 
         \begin{align*}
         & d_\alpha (x, y ) < \epsilon \\
         \implies & \inf \{ t > 0 : P( x, y, t ) < \alpha \} < \epsilon \\
         \implies &  P( x, y, \epsilon ) < \alpha \\
         \implies & y \in B ( x , \alpha, \epsilon ). 
         \end{align*}
         Hence $ B_\alpha (x, \epsilon ) \subset B ( x , \alpha, \epsilon )  $.  \\
         For the converse part, let $ G \in \tau_P  $. Then for every $ x \in G, $ there exists $ \alpha > 0 $ and $ t_0 > 0 $ such that 
         $ B(x, \alpha, t_0 ) \subset G $. Now, 
         \begin{align*}
         & y \in B(x, \alpha, t_0 ) \\
         \implies & P ( x, y, t_0 ) < \alpha \\
         \implies & P ( x, y, t' ) < \alpha, ~ \forall ~ t' > t_0 \\ 
         \implies & \inf \{ s > 0 : P( x, y, s ) < \alpha \} \leq t_0 \\
          \implies & d_\alpha (x, y ) \leq t_0
       \end{align*} 
    Choose $ s_0 > 0 $ such that $ s_0 > t_0 $. Then,      
       $  d_\alpha (x, y ) < s_0 $ and hence $  y \in B_\alpha (x, s_0 ) $. \\
         Therefore, 
           $ B(x, \alpha, t_0 ) \subset   B_\alpha (x, s_0 ) $ where $ s_0 > t_0 $. \\
           This completes the proof.
        \end{proof}
\end{thm}


\textbf{Conclusion:} In this article, we  exercise on topological properties, completeness etc. of generalized parametric metric space. All the results with respect to the general binary operation `o' are established. Moreover, we bring out that when the binary operation `o' is taken as `max',  both the topologies $ \tau_{d_\alpha} $ and $ \tau_P $ are identical and hence for those generalized parametric metric spaces topological results are same as metric space.  We hope the results of this article can significantly contribute  for further development in generalized parametric metric space.  
	~\\
	~\\
\textbf{Acknowledgment:}  
   The author AD  is grateful to University Grant Commission (UGC), New Delhi, India for awarding  senior research 
  fellowship [Grant No.1221/(CSIRNETJUNE2019)]. The authors are also thankful to the Department of Mathematics, Siksha-Bhavana, Visva-Bharati.  
	~\\
	~\\

%

\end{document}